\documentclass[a4paper,11pt]{article}
\usepackage[hypertexnames=false]{hyperref}
\usepackage{amsthm}
\usepackage{authblk}
\usepackage[letterpaper,margin=2cm]{geometry}
\usepackage{enumerate}
\usepackage{colonequals}

\usepackage{tikz,xcolor,amsfonts,amsmath,amssymb,mathrsfs,setspace}
\usepackage{hhline}

\newtheorem{thm}{Theorem}[section]
\newtheorem{dfn}[thm]{Definition}
\newtheorem{pro}[thm]{Proposition}
\newtheorem{cor}[thm]{Corollary}
\newtheorem{lemme}[thm]{Lemma}
\newtheorem{rmk}[thm]{Remark}

\newcommand{\di}{\operatorname{div}}

\newcommand{ \p }{\subset  }
\newcommand{\bd}{\begin{dfn}}
\newcommand{\ed}{\end{dfn}}
\newcommand{ \be}{\begin{equation}}
\newcommand{\ee}{\end{equation}}
\newcommand{\intom}{\int_\Omega }

\newcommand{\bt}{\begin{thm}}
\newcommand{ \et}{\end{thm}}
\newcommand{ \bp}{\begin{pro}}
\newcommand{ \ep}{\end{pro}}
\newcommand{ \bc}{\begin{cor}}
\newcommand{ \ec}{\end{cor}}
\newcommand{ \bl}{\begin{lemme}}
\newcommand{ \El}{\end{lemme}}
\newcommand{ \br}{\begin{rmk}}
\newcommand{ \Er}{\end{rmk}}

\newcommand{\s}{^{-1}}
\newcommand{ \R}{\mathbb{R}  }

\newcommand{\g}[1]{\mathbf #1}
\newcommand{\test}[1]{C_c^{\infty}(#1)}

\newcommand \Lim[2]{ \lim _{#1\rightarrow #2} }

\newcommand{\I}{\infty}
\newcommand{\el}[2]{  L^{#1}(#2)  }

\newcommand{\pd}{\partial}
\newcommand{\pair}[2]{\left\langle{#1},{#2}\right\rangle}

\newcommand{\reed}[1]{{\color{black}#1}}
\newcommand{\blu}[1]{{\color{black}#1}}

\newcommand{\pur}[1]{{\color{black}#1}}

\title{Uniqueness of Minimizers for Some Variational Problems Arising in Image Processing}
\author[1]{Romeo Awi}
\affil[1]{Department of Mathematics, Hampton University, Hampton, Virginia 23668\thanks{\textit{E-mail address}: romeo.awi@hamptonu.edu}}
\author[2]{Rohit Gupta}
\affil[2]{Department of Mathematics, Michigan State University, East Lansing, MI 48824\thanks{\textit{E-mail address}: guptaro1@msu.edu}}

\date{}

\begin{document}

\maketitle

\begin{abstract}
We will study an open problem pertaining to the uniqueness of  {minimizers for} a class of  {variational} problems emanating from Meyer's model  {for} the decomposition of an image into a geometric part and a texture part. 
Mainly, we are interested in the uniqueness of minimizers for the problem: 
\be \label{eq:abs-1}
 \inf \left\{J(u)+J^*\left(\frac{v}{\mu}\right): (u,v)\in
L^2({\mathrm{\Omega}})\times L^2({\mathrm{\Omega}}){,}\; f=u+v \right\}
 \ee
where the image $f$ is a square integrable function
on the domain ${\mathrm{\Omega}}$, the number $\mu$ is a parameter, the functional $J$ stands for 
the total variation and 
the functional $J^*$ is its Legendre transform.
We will consider Problem \eqref{eq:abs-1} as a special case of the problem: 
\be \label{eq:abs-2}
\inf\{s(f-u)+s^*(u): u\in \mathcal{X}\}
\ee
where $\mathcal X$ is a Hilbert space
containing $f$ and $s$ is a continuous semi-norm on $\mathcal X$.
In finite dimensions, we will prove that Problem \eqref{eq:abs-2} reduces to a projection problem
onto the polar of the unit ball associated to a given norm on an appropriate Euclidean space. 
We will also provide a characterization for the uniqueness of minimizers 
for a more general projection problem defined by using any norm and any nonempty, closed, bounded and convex set 
of an Euclidean space. Finally, we will provide numerical evidence in favor of the uniqueness of 
minimizers
for the decomposition problem. 
\end{abstract}

\section{Introduction}
An important problem in image processing consists of
decomposing a given image $ f\in \el2 {\mathrm{\Omega}} $ into a sum of a regular geometric part $ u $ and an oscillating texture part $ v $.
Most models in the literature
 use the total variation
\be
\el2{\mathrm{\Omega}} \ni w\mapsto J(w)=\intom 
 {\|D w\|}
=\sup\left\{
\intom w \di\phi dx: \phi\in \test {\mathrm{\Omega}},\;  {\|\phi\|_{\el \I {\mathrm{\Omega}}}}\leq 1
\right\}\nonumber
\ee
as well as its polar
\be\el2{\mathrm{\Omega}} \ni w\mapsto J^\circ(w)=\sup\left\{
\intom w g dx: g\in L^2({\mathrm{\Omega}}){,}\; J(g)\leq 1
\right\}.\nonumber
\ee
For instance, 
Rudin  {et al.} proposed the following restoration model in \cite{rof}:
\be\label{eq:rof model}
\inf\left\{\frac{1}{2\lambda}\|v\|^2 _{\el 2 {\mathrm{\Omega}}}+\intom  {\|D u\|}
:(u,v)\in  \el 2 {\mathrm{\Omega}}\times \el 2 {\mathrm{\Omega}}{,}\; f=u+v
\right\}.
\ee
The study of the model in \eqref{eq:rof model} led Meyer to propose in \cite{meyer} the model:
\be\label{eq:meyer model}
\inf\left\{J(u)+\alpha J^\circ (v):(u,v)\in \el 2 {\mathrm{\Omega}}\times \el 2 {\mathrm{\Omega}}{,}\; f=u+v
\right\}.
\ee
In \eqref{eq:meyer model}, 
\reed{we have a weighting parameter}
\blu{$\alpha$} 
 {and $J^\circ$ is the polar of $J$}.
Meyer's model in turn is approximated in \cite{auj} and \cite{vo} by using the following functional depending on the parameters $ \lambda $ and $ \mu $:
\be
F_{\lambda,\mu}: L^2({\mathrm{\Omega}})\times L^2( {\mathrm{\Omega}})\ni (u,v)\mapsto 
J(u)+\frac{1}{2\lambda}\| f-u-v\|^2 _{\el 2 {\mathrm{\Omega}}}
+\chi _{\mu}(v)\nonumber
\ee
where the functional $ \chi_\mu $ is defined over $ \el 2 {\mathrm{\Omega}} $ by
\be
\chi _{\mu}(v)=\begin{cases}
	0, &\text{if}~J^\circ (v)\leq\mu\\
  {\I}, &\text{if}~J(v)>\mu.
\end{cases}\nonumber
\ee
 {The approach taken by} \blu{Aujol et al{.}} in \cite{auj} consists of  {studying}  {the} problem:
\be\label{eq:aujol}
\inf\{F _{\lambda,\mu}(u,v): (u,v)\in
L^2({\mathrm{\Omega}})\times L^2( {\mathrm{\Omega}})\}.
\ee
In this paper, we are interested  {in}  {studying} the asymptotic case of Problem  
\eqref{eq:aujol}{, i.e.,} when $ \lambda $ tends to zero.
 We are then concerned
 with the problem: 
 \be \label{eq:main intro-1}
 \inf \left\{J(u)+J^*\left(\frac{v}{\mu}\right): (u,v)\in
L^2({\mathrm{\Omega}})\times L^2( {\mathrm{\Omega}}){,}\; f=u+v \right\}
 \ee
with $ J^* $ being  the Legendre transform of $ J $. 
 {Determining the u}niqueness 
 {of minimizers} in  {P}roblem \eqref{eq:main intro-1} is open  {and also} in  {its corresponding} discrete version:
\be \label{eq:main discrete intro-1}
 \inf \left\{J_d(u)+J_d^*\left(\frac{v}{\mu}\right): (u,v)\in
X\times X{,}\; f=u+v \right\}
 \ee
where $ X $ is the set of  {$ N\times N $}  {real} matrices and $ J_d $
is 
 {a}
discretization of $ J $ that we will define  {shortly}.
First we will define a discrete gradient operator
$ \nabla:X\to X\times X;\; u\mapsto \nabla u $   defined 
by 
$ (\nabla u) _{i,j}=((\nabla u) _{i,j}^1,(\nabla u) _{i,j}^2))
 $
 {for $i,j\in \{1,\dots,N\}$}  {and} \blu{where}
 \be
 (\nabla u) _{i,j}^1=\begin{cases}
	u _{i+1,j}-u _{i,j}
, &\text{if}~{1\leq} i<N\\
0
, &\text{if}~i=N\\	
\end{cases}\nonumber
\ee
and
\be
 (\nabla u) _{i,j}^2=\begin{cases}
	u _{i,j+1}-u _{i,j}
, &\text{if}~{1\leq}j<N\\
0
, &\text{if}~j=N.\\	
\end{cases}\nonumber
\ee
Next, for $ u\in X $, its discrete total  variation is defined by
 $ J_d(u)=\sum _{i,j=1}^ {N}  {\|}( \nabla u) _{i,j} {\|}.
 $
We refer the reader to \cite{aub} for more details on the decomposition problem. 
While  {the} existence 
 {of a minimizer} in Problems 
\eqref{eq:main  intro-1}  and
\eqref{eq:main discrete intro-1} 
follows directly from standard arguments in the calculus of variations,
 {proving} \reed{the} uniqueness is a challenge  {that} 
we will address in this paper.

Note that as far as  {the} uniqueness  {of minimizers} is concerned, the study
of Problems \eqref{eq:main  intro-1}  and
\eqref{eq:main discrete intro-1} may be reduced{, respectively,} to the study
of  {the following problems:}
\be \label{eq:main intro}
 \inf \left\{J(u)+J^*\left(v\right): (u,v)\in
L^2({\mathrm{\Omega}})\times L^2( {\mathrm{\Omega}}){,}\; f=u+v \right\}{,}
 \ee
\be \label{eq:main discrete intro}
 \inf \left\{J_d(u)+J_d^*\left(v\right): (u,v)\in
X\times X{,}\; f=u+v \right\}
 \ee
and this could in turn bring valuable
tools in the calculus of variations that may be helpful in  {determining} the uniqueness 
 {of minimizers} in 
problems involving the total variation
(see\blu{,} for instance\blu{,} \cite{awi-sed}){.}
We will first consider
Problem 
\eqref{eq:main  intro}
as a special case
of the more general problem:
\be\label{eq:intro *}
\inf\{s(f-x)+s^*(x):x\in \mathcal{X}\}
\ee
where $ \mathcal{X} $ is  {a Hilbert} space containing $ f $, the map $ s $ is a semi-norm 
on $ \mathcal{X} $
 {and} for all $ x\in \mathcal{X} $,
$ s^*(x)=\sup\{\pair xy-s(y): y\in \mathcal{X}\} $.
Let $ \mathcal{Y} $ be the orthogonal
complement of the  {subspace}
$ \{x\in \mathcal{X}:s(x)=0\} $, let $ \rho $ be the norm obtained
by restricting $ s $ to $ \mathcal{Y} $.
We will show in Lemma \ref{lem:semi-norm to norm}
that for every $ f\in \mathcal{X} $, we 
can find 
 {some} 
$ f_0\in \mathcal{Y} $ \blu{such that} 
Problem \eqref{eq:intro *}
is equivalent to  {the problem:}
\be
\inf\{\rho(f_0-x)+\rho^*(x):x\in \mathcal{Y}\}\nonumber
\ee
with
$ \rho^*(x)=\sup\{\pair xy-s(y): y\in \mathcal{Y}\} $. Note that in the definition of $ \rho^* $, the supremum is taken over $ \mathcal{Y} $ while for $ s^* $, the supremum is taken over $ \mathcal{X} $.
We  {can then} proceed to study the problem:
\be\label{eq:intro **}
\inf\{\rho(f_0-x)+\rho^*(x):x\in \mathcal{Y}\}
\ee
where $ \mathcal{Y} $ is  {a Hilbert} space containing $ f $ and
$ \rho $ is a norm on $ \mathcal{Y} $.
We may see Problem \eqref{eq:intro **}
as
\be\label{eq:intro ***}
\inf\{\rho(f-x): x\in D\}
\ee
where $ D=\{x\in \mathcal{Y}: \pair xy \leq 1\; \text{whenever}\; y\in \mathcal{Y}\; \text{and}\; \rho(y)\leq 1\}$ is the polar
of the unit ball $ \pur{\bar B_\rho(0,1)}=\{x\in \mathcal{Y}:\rho(x)\leq 1\} $.
We note that Problem \eqref{eq:intro ***}
is a projection problem onto the
polar of the unit ball $\bar B_\rho(0,1)$.

{Proving the uniqueness of minimizers in} Problem \eqref{eq:intro ***} is far from being trivial. Obviously, if the unit ball of 
$ \rho $ is strictly convex,  {then the minimizer is unique but this condition is not necessary to guarantee uniqueness}. 
In fact if we consider $ \mathcal{Y}=\R^N $ and we let 
$ \rho $ be the $ l^1 $-norm defined for $ x\in \R^N $ by $ \|x\|_1=\sum _{i=1}^N |x_i|
 $, then $ D=\{x\in \R^N: \|x\|_\I=\sup\{|x_i|: i=1,\dots , N\}\leq 1\} $. In this case, 
Problem \eqref{eq:intro ***}
admit{s} a unique solution 
\reed{despite the fact that the unit ball of the  $l^1$-norm is not 
strictly convex}
(see  {L}emma \ref{lem:uniqueness case}).

We will {also show} that Problem
\eqref{eq:intro ***}
may have several minimizers.
This is the case\blu{,} for 
instance\blu{,} 
when
$ \mathcal{Y}=\R^2 $
and the norm is defined
by 
$ 
\rho(x,y)=|x|+\max(|x|,|y|) 
$  {(see Proposition \ref{pro:cex}).}

We will reveal a strong connection
between
the uniqueness 
 {of minimizers}
in  {P}roblem \eqref{eq:intro ***}
and the existence in the unit ball of $ \rho $ of an edge
that is orthogonal to a vert{ex}.
We will show in Proposition \ref{pro:fundamental}
that 
the existence in the unit ball of $ \rho $ of an edge
that is orthogonal to a vert{ex}
is a sufficient condition
for nonuniqueness.
In dimension two, this
condition turns out to be
also necessary (see Theorem \ref{thm:1}).
We conjecture that this result
can be extended
to  {any} dimension {greater than two}.

 {Going}  {b}ack to  {P}roblem \eqref{eq:main  intro},
the  {subspace} $ \{u\in \el 2 {\mathrm{\Omega}}: J(u)=0\}$
is the set of  {all} constant functions
and its orthogonal
complement 
is the set
$ G=\{u\in \el 2{\mathrm{\Omega}}:\intom u dx=0\} $. Call $ T $ the
norm obtained by restricting $ J $ to $ G $.
To finish, we will provide numerical
 {evidence} that  {supports our conjecture that} the minimizer
in problem 
\eqref{eq:main discrete intro}  {is unique}.
 {The numerical} experiments consist
of taking a point
$ f $,  {choosing} several random starting points and  {then using} a projected \blu{subgradient} algorithm
to see if \blu{\pur{the iterates generated by} the algorithm \pur{converge} to the same solution or not.}

\section{Definitions and Notations}
Before proceeding further, we will recall some definitions and also fix some notations:
\bd
Let $ {\mathrm{\Omega}}\p \R^N $ be an open set. Suppose $ u\in L^ 1_{\mathrm{loc}}( {\mathrm{\Omega}} )$. The total variation of $ u $ is  {given by}
\be
\intom  {\|Du\|}
=\sup\left\{
\intom  {u \di\phi dx}: \phi\in \test {\mathrm{\Omega}},\;  {\|\phi\|_{\el \I {\mathrm{\Omega}}}}\leq 1
\right\}\nonumber
\ee
and we will set $ \|u\| _{BV({\mathrm{\Omega}})}=\intom  {\|Du\|}
+\|u\| _{\el 1 {\mathrm{\Omega}}}
 $.
\ed
We refer the reader to \cite{mtfpf} for more detail{s} on the space $ BV({\mathrm{\Omega}}) $. 
Let us point out that if $ u:\R\to \R $
is  a {piecewise constant function} which has  {finitely many} jumps at $ x_1<x_2<\cdots <x_k $, then 
\be
 \|u\| _{BV({\mathrm{\Omega}})}=\sum _{i=1}^k|u(x_i+)-u(x_i-)|
+\|u\| _{\el 1 {\mathrm{\Omega}}}\nonumber
\ee
with $ u(x_i-)=\Lim x {x_i^-} u(x) $ and 
$ u(x_i+)=\Lim x {x_i^+} u(x) $.\\

{\noindent\bf{Definitions and Notations from Convex Analysis:}}
\begin{enumerate}[(i)]
\item  {In what follows}, the set $ E $ stands for \blu{a Euclidean space}.
\item Let $ \rho $ be a norm on \blu{$ E $, then}  {the associated \blu{open} unit ball will be denoted by $ B_\rho(0,1) $}.
\item Let $ {\mathrm{\Omega}}\p E $ be  convex  and symmetric with respect to the origin (i.e., if $ x\in {\mathrm{\Omega}} $, then $ -x\in {\mathrm{\Omega}} $). Suppose also that $ {\mathrm{\Omega}} $ has a nonempty interior. To $ {\mathrm{\Omega}} $ we will associate
the functional
\be
\rho_{\mathrm{\Omega}}(x)=\inf\{t>0: t\s x\in {\mathrm{\Omega}}\}.\nonumber
\ee
The functional $ \rho_{\mathrm{\Omega}} $ is called the  {Minkowski functional} or  {the gauge}
of $ {\mathrm{\Omega}} $. The functional $ \rho_{\mathrm{\Omega}} $ is a norm and its unit ball is
\be
B_{\rho_{\mathrm{\Omega}}}(0,1)={\mathrm{\Omega}}=\{x\in {\mathrm{\Omega}}: \reed{\rho_{\mathrm{\Omega}}(x)< 1}\}.\nonumber
\ee
\item To a norm $ \rho $ on $ E $, we will associate the  {dual norm} $\rho^\circ :E\to [0,\infty)$ defined by
\be
\rho^\circ (y)=\max\{x\cdot y: x\in E{,}\; \rho(x)\leq 1\}= {\inf} \{\lambda>0: x\cdot y\leq \lambda \rho(x)\;  {\text{for all}} \; x\in E\} {.}\nonumber
\ee
Clearly for all $ x,y\in E $, we have $ x\cdot y\leq \rho(x)\rho^\circ(y) $  {(s}ee \cite{Dac_dmcv} for further properties of  {the gauge}
and its polar).
\item The polar of a convex set $ {\mathrm{\Omega}}\p E $ is defined by
\be
{\mathrm{\Omega}}^\circ=\{y:x\cdot y\leq 1\; \text{for all} \;x\in {\mathrm{\Omega}}\}.\nonumber
\ee 
\item A point $ x\in {\mathrm{\Omega}}\p E $ is said to be an extreme point of $ {\mathrm{\Omega}} $ if $ 2x=x_1+x_2 $ for $ x_1,x_2\in {\mathrm{\Omega}} $ implies  {that} $ x=x_1=x_2 $. We say that $ {\mathrm{\Omega}}\p E $ is strictly convex when it is convex and every point  {on} the boundary is an extreme point of $ {\mathrm{\Omega}} $.
%
\item
Let $ T:E\to \bar \R $ be a function. The Legendre transform of $ T $
is $ T^*:E\to \bar \R $
defined by
\be
T^*(y)=\sup\{x\cdot y-T(x): x\in E\}.\nonumber
\ee
\item The characteristic function of the set $ {\mathrm{\Omega}}\p E $ is  {defined} by 
\begin{align*}
	\chi_{\mathrm{\Omega}}(x)=\begin{cases}
	0, &\text{if}~x\in {\mathrm{\Omega}}\\
  {\I}, &\text{if}~x\not\in {\mathrm{\Omega}}.
	\end{cases}
\end{align*}
%
\end{enumerate}

\br
We have $ \rho^*=\chi _{\bar B _{\rho^\circ}(0,1)}
 $, i.e., for all $ y\in E $
 \begin{align*}
	\rho^*(y)=\begin{cases}
	0, &\text{if}~\rho^\circ(x)\leq 1\\
 {\I}, &\text{if}~\rho^\circ(x)> 1.
	\end{cases}
\end{align*}
\Er

\noindent  {The results in this subsection can be \blu{found,} for instance\blu{,} in \cite{brezis2010functional}, 
\cite{Dac_dmcv}, \cite{Eke_cavp} and \cite{Roc_ca}}. 
We will finish this  {subsection} with the following definition.

\bd

A subset $D$ of a normed vector space $\mathcal X$ %
is said to be a proximinal (respectively, a Cheby{s}hev) set with respect to the norm  {$\rho$}
if for every $ x_0\in \mathcal X$  the problem: 
\be
\inf\{{\rho}(x-x_0): x\in D\}\nonumber
\ee
admits a  solution (respectively, a unique solution).

\ed

{\noindent\bf{Image Modelization:}}
\begin{enumerate}[(i)]
\item We will denote  {by} $ X ${,} the space $ \R ^{{N\times N}} $ of all  {$ N\times N $}  {real} matrices. We will endow $ X $ with 
	the scalar product
$
	 \pair u v _X=\sum ^{N} _{i,j=1}u _{i,j} v _{i,j}
$
and the norm $ \|u\|_X=\sqrt {\pair u u _X} $.
\item The space $Y  $ is defined to be $ X\times X $.
For $ g=(g^1,g^2)\in Y $, we will define
\be  {\|g\|_\I=\max\left\{\sqrt{(g^1 _{i,j})^2+(g^2 _{i,j})^2}: i,j=1,\dots,N\right\}}\nonumber.\ee
\item The  {discrete gradient} operator $ \nabla:X\to Y;\; u\mapsto \nabla u $ is defined  {for $i,j=1,\dots,N$} by 
$ (\nabla u) _{i,j}=((\nabla u) _{i,j}^1,(\nabla u) _{i,j}^2))
 $ 
 where
 \be
 (\nabla u) _{i,j}^1=\begin{cases}
	u _{i+1,j}-u _{i,j}
, &\text{if}~i<N\\
0
, &\text{if}~i=N\\	
\end{cases}\nonumber
\ee
and 
\be
(\nabla u) _{i,j}^2=\begin{cases}
	u _{i,j+1}-u _{i,j}
, &\text{if}~j<N\\
0
, &\text{if}~j=N {.}\\	
\end{cases}\nonumber
\ee
\item For $ u\in X $, its total  variation is defined by  $ J_d(u)=\sum _{i,j=1}^N  {\|}( \nabla u) _{i,j} {\|}.
 $
\item The  {discrete} divergence operator  {$\di: Y\to X$} is defined by
 the relation
 \be
 \pair {-\di (p)}{u}_X=\pair{p}{\nabla u}_Y\;  {\text{for all}} \;u\in X,\; p\in Y.\nonumber
 \ee
 {Note} that: 
 \be
(\di p) _{i,j}=(\di p) _{i,j}^1+(\di p) _{i,j}^2\nonumber
\ee
where
\be
(\di p)^1 _{i,j}=\begin{cases}
p^1 _{1,j}{,}
 &\text{if}~i=1\\
	p^1 _{i,j}-p^1_{i-1,j}{,}
 &\text{if}~1<i< {N}\\	
-p^1_{d-1,j}{,}
 &\text{if}~i= {N}\\
\end{cases}\nonumber
\ee
and
\be
(\di p)^2 _{i,j}=\begin{cases}
p^2 _{i,1}{,}
 &\text{if}~j=1\\
	p^2 _{i,j}-p^2_{i,j-1}{,}
 &\text{if}~1< {j}< {N}\\	
-p^2_{i,d-1}{,}
 &\text{if}~{j}= {N.}\\
\end{cases}\nonumber
\ee
\end{enumerate}

%
%
%
%
%
%
%
%
%
%
%
%
%
%
%

\section{{Minimizatio{n}}  {Problems}  {I}nvolving a  {S}emi-{N}orm and its Legendre  {{T}ransfor{m}}}
In this section, we will consider
two minimization problems,
one involving a semi-norm and the other a norm.
We will show that these problems
are equivalent or in other words, they have the same set of minimizers.

\bl\label{lem:semi-norm to norm properties}
Let $\mathcal  X $ be  {a Hilbert} space with  {the} scalar product $ \phi:\mathcal X\times\mathcal  X\to \R $. Let $ \rho:\mathcal X\to [0,\I) $ be a \reed{continuous} semi-norm. Define the set 
$ G $ by
\be
G=\{x\in\mathcal  X: \phi(x,y)=0\; \text{whenever} \;y\in\mathcal  X\; \text{and} \;\rho(y)=0\}.\nonumber
\ee
Call $ \rho_G $ the norm obtained by restricting $ \rho $ to $ G $ and let $ \rho_G^* $ be its Legendre transform,  {i.e.}, for all $ y\in G $ we have
\be 
{\rho_G^*(y)=\sup\{\phi(w,y)-\rho_G(w): w\in G\}}.\nonumber
\ee
 {Let $K_\rho =\{x\in\mathcal  X:\rho(x)=0\}$, then the} following  {claims} \blu{hold}:
\begin{enumerate}[(i)]
	\item 
For all $ (x,y)\in\mathcal  X\times K_\rho $, we have
	$
\rho(x+y)=\rho( {x})$.

\item We have $\rho^*(x)=\I$ whenever $x\not\in G.$

\item If $ x\in G ${,} then 
$\rho_G^*(x)=\rho^*(x)$.
\end{enumerate}
\El

\begin{proof}
 {Firstly note} that $ G $ is the orthogonal complement of $ K_\rho $:
\be
G=K_\rho ^{\perp}=\{x\in\mathcal  X: \phi(x,y)=0\; \text{for all} \;y\in K_\rho\}.\nonumber
\ee
For all $ x\in\mathcal  X $ we will denote  {by} $ \hat x ${,} the orthogonal projection of $ x $ onto the  {subspace} $ K_\rho $. We  {then} have
\be\label{eq:rho hat}
 \rho(\hat x)=0\;  {\text{for all}} \;x\in\mathcal  X .
 \ee
\begin{enumerate}[(i)]
\item Let $ x\in\mathcal  X $ and $ y\in K_\rho $. Using the fact that $ \rho $ is subadditive and relation \eqref{eq:rho hat}, we have
\be
\rho(x+y)\leq \rho(x)+\rho(y)=\rho(x)\leq \rho(x+y)+\rho(-y)=\rho(x+y){,}\nonumber
\ee
 {from which w}e deduce \blu{that:}
\be\label{eq:rho K}
\rho(x+y)=\rho(y)\;  {\text{for all}} \;x\in \mathcal X,\; y\in K_\rho.
\ee

\item Suppose $ x\not\in G $. It follows that $ \hat x\neq 0 $  {and w}e  {now} use the definition of $ \rho^* $  {and} \eqref{eq:rho hat} to  {obtain} for all $ t\in \R $:
\be
\rho^*(x)\geq \phi( {t\hat x,x}) {-}\rho(t\hat x)=t \phi(\hat x,\hat x).\nonumber
\ee
We let $ t $ go to  {infinity} to  {obtain} $ \rho^*(x)=\I $.
Whence,
\be\rho^*(x)=\I\; \text{whenever} \;x\not\in G.\nonumber\ee

\item Suppose $ x\in G $. We have
 {\begin{align*}
	\rho^*(x)&=\sup\{\phi(x,y)-\rho(y): y\in \mathcal X
	\}
	\\
&=\sup\{\phi(x,y-\hat y)-\rho(y): y\in  \mathcal X
	\}\; (\text{as $x\in G$})
	\\
&=\sup\{\phi(x,y-\hat y)-\rho(y-\hat y): y\in \mathcal X
	\}\; (\text{by \eqref{eq:rho K}})
	\\
&=\sup\{\phi(x,w)-\rho(w): w\in G
	\}
	\\
&=\rho_G^*(x).
\end{align*}}
\end{enumerate}
\end{proof}

%
\bl\label{lem:semi-norm to norm}
Let $\mathcal  X $ be  {a Hilbert} space with scalar product $ \phi:\mathcal X\times \mathcal X\to \R $. Let $ \rho:\mathcal X\to [0,\I) $ be a \blu{continuous} semi-norm. Define the set 
$ G $ by
\be
G=\{x\in \mathcal X: \phi(x,y)=0\; \text{whenever} \;y\in \mathcal X\; \text{and} \;\rho(y)=0\}.\nonumber
\ee
Call $ \rho_G $ the norm obtained by restricting $ \rho $ to $ G $ and let $ \rho_G^* $ be its Legendre transform,  {i.e.}, for all $ y\in G $ we have
\be
{\rho_G^*(y)=\sup\{\phi(w,y)-\rho_G(w): w\in G\}}.\nonumber
\ee
For every $ f\in \mathcal X $, there exists \blu{some} $ f_0\in G $ \blu{such that} the following problems are equivalent (i.e., they have the same set of minimizers):
\be\label{eq:min semi} \inf\{\rho(f-x)+\rho^*(x): x\in \mathcal X\}{,} \ee
\be\label{eq:min norm}
\inf\{\rho_G(f_0-x)+\rho_G^*(x): x\in G\}{.}\ee
\El

\begin{proof}
We will use  {the second claim in} Lemma \ref{lem:semi-norm to norm properties}  {t}o deduce that  {P}roblem
\eqref{eq:min semi} is equivalent to  {the problem:}
\be
\label{eq:dum 1}
\inf\{\rho(f-x)+\rho^*(x): x\in G\}.
\ee
 {Let}
$ K_\rho =\{x\in \mathcal X:\rho(x)=0\}$  and call $ \hat f $
the orthogonal projection 
of $ f $ on $ K_\rho $. We will choose $ f_0=f-\hat f $.
 {Using the first claim in} Lemma \ref{lem:semi-norm to norm properties},  {it follows that}  {P}roblem \eqref{eq:dum 1} is equivalent
to  {the problem:}
\be
\label{eq:dum 2}
\inf\{\rho(f_0-x)+\rho^*(x): x\in G\} .
\ee
Finally, we will use the fact that  for all $ x\in G $
 {we have} $ f_0-x \in G$ and 
 {the third claim in} Lemma \ref{lem:semi-norm to norm properties}
to deduce  {that}  {P}roblem \eqref{eq:dum 2} is equivalent
to  {P}roblem \eqref{eq:min norm}. We have  {now} established that  {P}roblem \eqref{eq:min semi} is equivalent
to  {P}roblem \eqref{eq:min norm}.
\end{proof}

\br
 {Let us make the following observations:}
\begin{enumerate}[(i)]
 \item 
We do not have uniqueness of $ f_0 $  in Lemma \ref{lem:semi-norm to norm}. An example  {of this} is
the case  {when} $ \blu{\mathcal X=G=\R^2}$  {and} $ \rho $ is the Euclidean
norm. If  $ f=(1,0) $,  {then} any $ f_0 $ of the form $ f_0=(t,0) $ with $ t\geq 1 $ will satisfy the conclusion of Lemma \ref{lem:semi-norm to norm}.
\item 
 {It is apparent from the proof of   Lemma \ref{lem:semi-norm to norm} that $f_0$ may be  {taken\blu{,}} for instance\blu{,} 
to be the orthogonal projection
of $f$ onto $G$.}

\end{enumerate}
\Er

\subsection{The  {C}ase of  {the} {T}otal  {V}ariation}
We will study  {the problem:}
 {\be
\inf\{J(f-v)+J^*(v): v\in \el2{\mathrm{\Omega}}\}.\nonumber
\ee}
Let $$
 G=\left\{v\in \el 2 {\mathrm{\Omega}}: \intom u dx=0\right\}\; \text{and} \;{ G_d=\left\{u\in X: \sum_{i,j=1}^N u_{i,j}=0\right\}}.
 $$
Call $ T:G\to [0,\I) $ the restriction of $ J $ to $ G $
 {and $ T_d:G_d\to [0,\I) $ the restriction of $ J_d $ to $ G_d $}.
 {The next Corollary follows from Lemma \ref{lem:semi-norm to norm}.}
\bc\label{cor:}
For every $ f\in \el 2{\mathrm{\Omega}} $, there exists   {some} $ f_0\in G $
\blu{such that} the following problems are equivalent  {(i.e., they have the same set of minimizers)}:
\be\inf\{J(f-x)+J^*(x): x\in \el 2{\mathrm{\Omega}}\}{,}\nonumber \ee
\be
\inf\{T(f_0-x)+T^*(x): x\in G\}.\nonumber\ee
\ec

\br
 {A result analogous to Corollary \ref{cor:} also holds for the discrete version, i.e., $\el 2{\mathrm{\Omega}}$, $G$ and $J$ replaced with $X$, $G_d$ and $J_d$, respectively}. 
\Er

\bp
Neither of the unit balls of the norms $ T $  and $T_d$	is strictly convex.
\ep

\begin{proof}
For $ (a,b)\in \R^2 $   define the matrix $\bar M(a,b) $ by $ [\bar M(a,b)]_{1,j}=a $ for $ j=1, {\dots},N $
and $[\bar M(a,b)]_{i,j}=b $ for $ i=2, {\dots},N $ and $ j=1, {\dots},N $.
Let $ \hat M(a,b) $ be the constant matrix with
coefficient $ Na+N(N-1)b $ and $ M(a,b)=\bar M(a,b)-\hat M(a,b) $.
It  {follows} that $ M(a,b)\in G {_{d}} $ and $ T {_{d}}(M(a,b))=N|b-a| $. \blu{One can now verify that the} following equations hold:
\begin{align*}
	2M\left(\frac12,-\frac12\right)&=M(1,0)+M(0,-1){,}
	\\
	2T {_{d}}\left(M\left(\frac12,-\frac12\right)\right)&=T {_{d}}\left(M(1,0)\right)+T {_{d}}\left(M(0,-1)\right).
\end{align*}
Hence the unit ball of $ T {_{d}} $ is not strictly convex.  {A similar proof also holds for the other case.}
\end{proof}

\bl\label{lem:red 2}

Suppose that:
\begin{enumerate}[(i)]
	\item The matrix $f_0=f-\hat f$  {where} $ \hat f $ is the $  {N\times N} $  {constant} matrix
with  {coefficient}
$  {\sum _{i,j = 1}^{N} f _{ij}}
 $. 
 \item The set 
 $ 
 D _{d}=\{y\in G_d: T _d^\circ(y)\leq 1\}=
 \reed{\bar B}_{T_d^\circ}(0,1)
 $.
\end{enumerate}
Then
Problem \eqref{eq:main discrete intro}
is equivalent to  {the problem:}
\be\label{eq:reduced}\inf\{T {_{d}}(f_0-v):v\in D {_{d}}\}.\ee
\El

\begin{proof}
As $ T _d $ is a norm on $ G_{d} $,
we have $ \blu{T_d^*}=\chi_{\reed{\bar B} _{T _{d}^\circ}(0,1)
} $. We will next use  {Corollary \ref{cor:}},
to deduce Lemma \ref{lem:red 2}.
\end{proof}

\br
Problem \eqref{eq:reduced} is  {a} projection  {problem} with respect to the  {norm} $ T_\reed{d} $  {o}nto the set $ D_\reed{d} $ which is the dual of the unit ball associated to $ T_\reed{d} $.
\Er

\section{The  {P}rojection {Problem} onto the  {D}ual  {U}nit  {B}all}
{Let} $ \rho $ be a norm on \blu{$E$}  {and} $ x_0\in E $.
We  {will} study the problem:
\be \label{eq:pp}
\inf\{\rho(x_0-x)+\rho^*(x): x\in E\}
\ee
where the map $\rho^*  $ is the Legendre 
 {transform} of $ \rho $ and is defined %
by
\be
\rho^*(y)=\sup\{x\cdot y -\rho(x): x\in E\}\; \text{for all} \; y\in E.\nonumber
\ee
Problem \eqref{eq:pp} is a projection problem with respect to the norm $\rho$ onto the dual of the unit ball associated to $\rho$. To see this, recall that for all $ x,y\in E$, one has
$ x\cdot y\leq \rho(x)\rho^\circ(y) $  {and it}
holds that 
 $ \rho^*=\chi _{{\reed{\bar B} _{\rho^\circ}(0,1)}}
 $.
Let $ D=\{x\in E: \rho^*(x)=0\} ${,}  {t}hen $ D=\{x\in E:\rho^\circ(x)\leq 1\}=\{x: x\cdot y \leq \rho(y)\;  {\text{for all}} \;y\in \R^N\} $ \blu{and it follows that:}
\be \label{eq:ppg}
\inf\{\rho(x_0-x)+\rho^*(x): x\in E\}=\inf\{\rho(x_0-x): x\in D\}.
\ee
Furthermore, as $D \p E $ is closed and bounded and $ \rho $ is a norm, we read from \eqref{eq:ppg} that    Problem \eqref{eq:pp} admits a minimizer.


\bl\label{lem:proj1}
 {The following claims hold:}
\begin{enumerate}[(i)]
	\item If $ x_0\in D $, \blu{then} Problem  \eqref{eq:pp} admits a unique solution which is $ x_0 $.
	\item If $ x_0\not\in D $, \blu{then} any minimizer of Problem  \eqref{eq:pp} lies on  {t}he boundary of $ D $.
	\item If $ x_0\not\in D $ and $ D $ is strictly convex (i.e., every point  {on the boundary of $D$} is an extreme point of $ D $){,} then   Problem  \eqref{eq:pp}  admits a unique minimizer.
\end{enumerate}
\El

\begin{proof}
We will prove the claims of the lemma one by one:
\begin{enumerate}[(i)]
\item Suppose $ x_0\in D $. As $ \rho(x_0-x)\geq 0=\rho(x_0-x_0) $ for all $ x\in E $, it holds that $ x_0 $ is a minimizer of Problem  \eqref{eq:pp}.
Let $ x_1 $ be another minimizer of Problem  \eqref{eq:pp}. Then $ \rho(x_0-x_1)=0 $, this implies that $ x_0=x_1 $. Hence the minimizer is unique.
\item Suppose $ x_0\not\in D $ and suppose  \pur{also} {that} a minimizer $ x_1\in D $ lies in the interior of $ D $. We may find $ \blu{\epsilon\in(0,1)}$ \blu{such that}
the  ball $ B_\rho(x_1,\epsilon)\p D $.  Since 
$ x_0\not\in D $ and $ x_1\in D $, we have $ \rho(x_0-x_1)>0 $.
Let \reed{$t\in(0,1)$ be such that} $ 0<t\rho(x_0-x_1)<\epsilon $ and set $ x_t=x_1+t(x_0-x_1) $. We have that $ x_t\in D $ and 
\be
 {\rho(}x_0-x_t)= {\rho(}x_0-x_1-t(x_0-x_1))= {\rho(}(1-t)(x_0-x_1))=(1-t) {\rho(}x_0-x_1).\nonumber
\ee
Thus \be\label{eq:interior} (1-t) {\rho(}x_0-x_1)\geq  {\rho(}x_0-x_1) .\ee
 {The} inequality \eqref{eq:interior} reads
$ 1-t\geq 1 $, which is  {absurd} as $ t>0 $  {and}  {h}ence $  {x_1\in \partial D}$.
\item Let $ x_1,x_2 $ be  {two distinct} minimizers of Problem \eqref{eq:pp}. Then $ x_1,x_2 $ lie on the boundary  {of $D$} and are extreme points. It holds that $\frac{x_1+x_2}2$ is also a
minimizer and must be an extreme point. Hence $ x_1=x_2 $ and the minimizer is unique.
\end{enumerate}
\end{proof}

%
%
%
%
%
%

\subsection{A  {C}ase  {where the Minimizer is Unique}:  {T}he  {$  {l^p(\R^N)} $-{N}orms}}
If $ p\in (1,\I) $,  {then}  {the} unit ball of the \pur{$ l^p(\R {^N}) $-norm} is strictly convex,
hence any nonempty, closed, bounded and convex set is a Cheby{s}hev set. In particular, the closed dual unit ball $\bar B_{\rho^\circ}
(0,1) $ is  {also}  {a} Cheby{s}hev  {set}.
 We will study  {next} the case of  the \pur{$ l^1(\R {^N}) $-norm}.
Consider for this purpose the map $ \alpha:\R\to \R $ defined by
\be
\alpha(t)=
\begin{cases}
	-1, &\text{if}~t\leq -1\\
	t, &\text{if}~-1<t<1\\
	1, &\text{if}~t\geq-1 {.}\\	
\end{cases}\nonumber
\ee
Call $ P:\R {^N}\to \R {^N} $,  the map defined for $ x\in \R {^N} $  by 
$ (P(x))_i=\alpha(x_i) $ for  $ i=1, {\dots}, {N} $. 
\bl\label{lem:uniqueness case}
Let $ f\in \R {^N} ${, then}  {t}he unique solution of the problem:
 {\be\label{eq:minL1}
 \inf\{\|f-u\|_1: u\in \R {^N},\; \|u\|_\I\leq 1\} 
\ee}
is given by $ P(f) $.
\El

\begin{proof}
 {Observe} that  for $ t\in \R $, $  {\alpha(t)} $ is the unique solution of the 
problem:
$$  {\inf\{|t-s|: s\in \R,\; | s|\leq 1\}.}
$$
For $ f,u\in  {\R^N} $ and $ \|u\|_\I\leq 1 $, one has
\be
	\|f-u\|_1=\sum_{i=1}^ {N}|f_i-u_i|
	\geq\sum_{i=1}^ {N}|f_i-\alpha(f_i)|\nonumber
\ee
with equality if and only if $ u_i=\alpha(f_i) $ for all $ i=1, {\dots},  N $.
We deduce that the unique minimizer of 
Problem \eqref{eq:minL1} is  {given by} $ P(f) $.
\end{proof}

%

%
%

%
%
%
\subsection{A  {C}ase  {where the Minimizer is not Unique}:  {A}  {N}orm for which the  {D}ual  {U}nit  {B}all is not Cheby{s}hev}
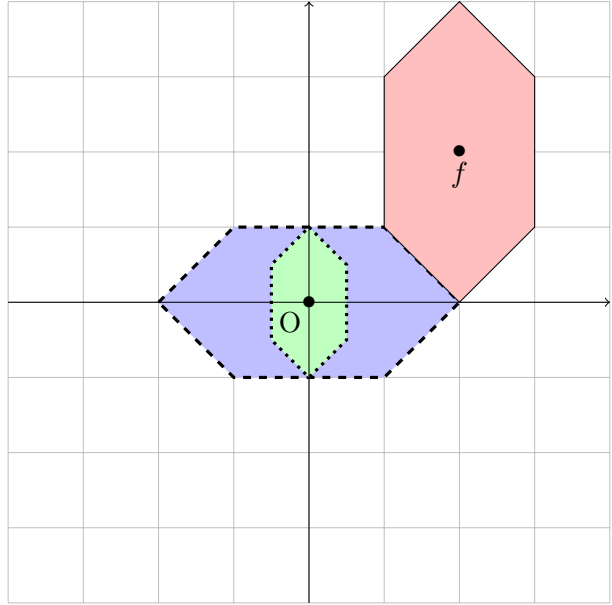
\begin{figure}[h]
\begin{center}
\begin{tikzpicture}
\draw[step=1cm,lightgray,very thin] (-4,-4) grid (4,4);
\draw[very thick,dashed,fill=blue!25] (-2,0)--(-1,1)--(1,1)--(2,0)--(1,-1)--(-1,-1)--(-2,0);
\draw [very thick,dotted,fill=green!25] (0,1)--(.5,.5)--(.5,-.5)--(0,-1)--(-.5,-.5)--(-.5,.5)--(0,1);
\draw[fill=red!25] (1,1)--(2,0)--(3,1)--(3,3)--(2,4)--(1,3)--(1,1);
\draw[->] (-4,0) -- (4,0);
\draw[->] (0,-4) -- (0,4);
\draw  (2,2) node {$\bullet$} node[below]{$f$};
\draw  (0,0) node {$\bullet$};
\draw  (-0.25,0) node[below]{O};
\end{tikzpicture}	
\end{center}
\caption{A norm with a non-Chebyshev dual. The polygons with dotted, dashed and solid boundaries are the unit ball, dual unit ball and the ball of radius two centered at $f=(2,2)$, respectively.}
\label{fig:1}
\end{figure}

 {Let} $ N=2 $  {and}  {c}onsider the polygon $ A \p \R {^N}$ with vertices  {located at} 
$(0,1)$, $(0.5,0.5)$, $(0.5,-0.5)$, $(0,-1)$, $(-0.5,-0.5)$ and $(-0.5,0.5)$. The polygon $ A $ is also characterized
by the  {following} inequalities:
\begin{align*}
	x+y\leq &1{,} & 2x\leq& 1{,} & x-y\leq& 1{,}\\
	-x+y\leq &1{,} & -2x\leq &1{,} & -x-y\leq &1.
\end{align*}
It follows that the  {polar} of $ A $ is the polygon $ A^\circ $
with vertices  {located at} 
$(1,1)$, $(2,0)$, $(1,-1)$, $(-1,-1)$, $(-2,0)$ and $(-1,1)$ ({s}ee  {F}igure \ref{fig:1}).
To $ A $ we will associate the norm 
$
\rho_A(\g x)=\inf\{t: t ^{-1}\g x\in A\}
$
which is the  {Minkowski functional or} gauge  {function} associated to $ A $. As a consequence 
\be
\rho^\circ(\g x)=\inf\{t: t^{-1}\g x\in A^\circ\}.\nonumber
\ee
 {Note} \blu{that:}
\be
\rho(x,y)=|x|+\max(|x|,|y|)=\begin{cases}
	2|x|, &\text{if}~|x|>|y|\\
	|x|+|y|, &\text{if}~|x|\leq|y|\\
	
\end{cases}\nonumber
\ee
while
\be
\rho^\circ(x,y)=\frac12|y|+\frac12\max(|x|,|y|)=\begin{cases}
	|y|, &\text{if}~|y|>|x|\\
	\frac12(|x|+|y|), &\text{if}~|y|\leq|x| {.}\\
	
\end{cases}\nonumber
\ee

\bl \label{lem:dub}
The dual unit ball $ A^\circ $ of $ \rho $ is not a Cheby{s}hev set.
\El

\begin{proof}
Consider $ f=(2,2) ${, then the minimizers of the problem:}
\be
\inf \{\rho(f-u): u\in A^\circ\}\nonumber
\ee
are {points of the form} $ (1+t,1-t) $  {with} $ t\in [0,1] $. Hence
$ A^\circ $ is not  {a} Cheby{s}hev  {set} with respect to the norm $ \rho $.
\end{proof}

\reed{
The following \blu{result} is a consequence of Lemma \ref{lem:dub}.
\bp \label{pro:cex}
\blu{Let $\rho$ be the norm on $\R^2$ defined by $\rho(x,y)=|x|+\max(|x|,|y|)$. There exists some $f_0\in \R^2 $ such that
the problem:
\be
\inf\{\rho(f_0-x)+\rho^*(x): x\in \R^2\}\nonumber
\ee
admits infinitely many solutions. In particular, the closed dual unit ball associated to $\rho$ 
is not a Chebyshev set.}
\ep
}

\subsection{Characterization for the Uniqueness of Minimizers}
We will first provide a characterization for the uniqueness of  minimizers for a projection problem defined by using any norm $\rho$ on $E$ and any nonempty, closed, bounded and convex set ${\mathrm{\Omega}}$ of $E$. In order to state subsequent results, we will have the need to fix some additional notations:
\begin{enumerate}[(i)]
\item{
For $ a\in E\setminus\{0\} $ and $ \alpha\in \R $, we will define the hyperplane:
\be
H(a,\alpha)=\{x\in E: x\cdot a=\alpha\}\nonumber
\ee
and the closed halfspaces:
\begin{align*}
H^+(a,\alpha)&=\{x\in E: x\cdot a\geq \alpha\},\\
H^-(a,\alpha)&=\{x\in E: x\cdot a\leq\alpha\}.
\end{align*}}
\item{The line segment between $a,b\in E$ is denoted by
\be
[a,b]=\{a+t(b-a): t\in [0,1]\}.\nonumber
\ee}
\item{{F}or a  {convex} set \blu{$ {\mathrm{\Omega}}\p E$}, we will denote by $  {\textnormal{extr}}({\mathrm{\Omega}}) $ the set of its extreme points.}
\item{The unit sphere associated with the norm $\rho$ on $E$ will be denoted by $  {\pd \bar B_\rho}(0,1) $.}
\end{enumerate}



%
\bp\label{pro:fundamental}
Let $ {\mathrm{\Omega}} $ be a \pur{nonempty, closed, bounded and convex subset of $ E $ and assume $\rho$ is a norm on $E$.} Then the following statement{s} are equivalent:
\begin{enumerate}[(i)]
	\item There exists  {some} $ x_0\in E $ for which 
	the problem: 
	\be
	\inf \{\rho(x_0-x): x\in {\mathrm{\Omega}}\}\nonumber
	\ee
	has more than one minimizer.
\item There \blu{exists} \reed{some vector} $ a\in E $, \blu{distinct points $ w_1,w_2\in{\mathrm{\Omega}} $, $ u_1,u_2\in {\pd \bar B_\rho}(0,1) $} and $ {r}\neq 0 $ such that 
the following three conditions hold:
\begin{enumerate}[(a)]
	\item We have ${\mathrm{\Omega}} \p H^+(a,w_1\cdot a) $
	 and $ [w_1,w_2] \p H(a,w_1\cdot a) $.
	 \item We have $ B_\rho(0,1) \p H^-(a,u_1\cdot a) $
	 and $ [u_1,u_2] \p H(a,u_1\cdot a) $.
	 \item We have $ w_1-w_2= {r}(u_1-u_2) $.
\end{enumerate}
\end{enumerate}
\ep

\begin{proof}
 {We will first prove that the first statement implies the second one.}
Suppose $ w_1, w_2 $ are two distinct minimizers of  {the problem:} 
\be
\inf\{\rho(x_0-x):x\in {\mathrm{\Omega}}\}\nonumber
\ee
and let $r\colonequals\inf\{\rho(x_0-x):x\in {\mathrm{\Omega}}\}>0$. Then $ [w_1,w_2] $ 
is a set of minimizers and 
$ [w_1,w_2] \p \pd {\mathrm{\Omega}}$.  {Since} $ {\mathrm{\Omega}} \cap B_\rho(x_0,r) =\emptyset$, we  {can then} use the theorem  {on the} separation of convex sets to
find $ a\in E $ and $  {\alpha\in\R} $ \blu{such that}
$ {\mathrm{\Omega}}\p H^+(a,\alpha) $
and $ B_\rho(x_0,r)\p H^-(a,\alpha) $. Next, as $  {\pd \bar B_\rho}(x_0,r)\p H^-(a,\alpha) $ and
$[ w_1,w_2]\p {\mathrm{\Omega}}\cap   {\pd \bar B_\rho}(x_0,r) $, we deduce that
$[ w_1,w_2]\p H(a,\alpha) $ and $  {\alpha=w_1\cdot a=w_2\cdot a} $.
Furthermore, let us consider the  {affine} map $ T:E\to E $
defined for $ x\in E $ by $ T(x)=r\s(x-x_0) $ and let $ u_1=T(w_1) $ and $ u_2=T(w_2) $.
For $ t\in [0,1] ${, we have}
\be
 {(}u_1+t(u_2-u_1) {)}\cdot a= {(}u_1+tr\s(w_1-w_2) {)}\cdot a=u_1\cdot a.\nonumber
\ee
Thus $ [u_1,u_2]\p H(a,a\cdot u_1) $  {and also} $ u_1,u_2\in \bar B_\rho(0,1) $. Let $ x\in \bar B_\rho(0,1) $ and set $ y=T\s(x) $. We have  $ y\in B_\rho(x_0,r) $ and then $ y\cdot a\leq  {w_1\cdot a} $.  {Since} $ y=rx+x_0 ${, we have}
\begin{align*}
	(rx+x_0)\cdot a&\leq  {w_1\cdot a}\\
	r x\cdot a&\leq  {(w_1-x_0)\cdot a}\\
	x\cdot a&\leq  {r^{-1}(w_1-x_0)\cdot a}\\
	            &= {u_1\cdot a}.
\end{align*}
Therefore, $ \bar B_\rho(0,1)\p H^-(a, {u_1\cdot a}) $.
Finally, $u_1-u_2=r\s(w_1-w_2)  $.

 {We will now prove that the second statement implies the first one.} 
Let $ x_0=w_1- {r}u_1 $  {and we will now} show that
$ w_1 $ and $ w_2 $ are   {two} distinct  minimizers
of  {the problem:}
\be
\inf\{\rho(x_0-x):x\in {\mathrm{\Omega}}\}.\nonumber
\ee
\blu{Firstly}  {note} that $ x_0=w_2- {r}u_2 $  {and}  {f}or $ i=1,2$, we have
\be
\rho(w_i-x_0)=\rho( {r}u_i)= {r}\rho(u_i)= {r}\nonumber
\ee
which implies $ w_1,w_2\in  {\pd \bar B_\rho}(x_0, {r}) $.
Suppose $ \rho(x-x_0)< {r} ${,}  {t}hen $  {r}\s(x-x_0)\in B_\rho(0,1) $. As $ B_\rho(0,1)\p H^-(a,u_1\cdot a) $, we \blu{have}
\begin{align*}
	 {r}\s(x-x_0)\cdot a&\leq u_1\cdot a
	\\
		 {r}\s(x-w_1+ {r}u_1)\cdot a&\leq u_1\cdot a
		\\
 {r}\s(x-w_1)\cdot a&\leq 0
\\
x\cdot a&\leq w_1\cdot a.
\end{align*}
Thus $ B_\rho(x_0, {r})\p H^-(a,w_1\cdot a) $
 {and}  {a}s $ w_1\in  {\pd \bar B_\rho}(x_0, {r}) $, we deduce that 
$ H(a,w_1\cdot a) $ is a supporting hyperplane
of $ B_\rho(x_0, {r}) $.
Thus $ x\in H^+(a,w_1\cdot  {a}) $ is equivalent to
$ \rho( {x-x_0})\geq  {r} $.
As $ {\mathrm{\Omega}}\p  H^+(a,w_1\cdot  {a})$, we deduce 
that $ w_1 $ and $ w_2 $ are  {two distinct} minimizers of  {the problem:}
\be
\inf\{\rho(x_0-x):x\in {\mathrm{\Omega}}\}\nonumber
\ee
{since}  {b}y assumption 
$ w_1 $ and $ w_2 $ are distinct.
\end{proof}

\bl\label{lem:nonuniqueness=>}
Let $ \rho $ be a norm on \blu{the $N$-dimensional Euclidean space} $  {E} $  {\blu{such that} its closed  unit ball
is the convex hull of a finite number of points.} 
If the  {p}roblem: 
\be
\inf \{\rho(x_0-x)+\rho^*(x):x\in E\}\nonumber
\ee
admits a unique  {minimizer} for all $ x_0\in E ${,}  then the following
set: 
\begin{gather}
W=\left\{(x_1,x_2,x_3)\in\textnormal{extr}(B _{\rho}(0,1))\times\textnormal{extr}(B _{\rho}(0,1))\times\textnormal{extr}(B _{\rho}(0,1)),\right.\nonumber\\ 
\left[x_2,x_3]\p  {\pd \bar B_\rho}(0,1),\; x_2\neq x_3,\; x_1\cdot(x_2-x_3)=0\right\}\nonumber
\end{gather}
is empty.
\El

\begin{proof}
Assume $(x_1,x_2,x_3) {\in W} $ and  {let} $ D $  {be the} polar of $ \reed{\bar B}_\rho(0,1) $.
As 
$ x_1 $ is an extreme point of $\reed{\bar B}_{{\rho}}(0,1) $, 
 {there exists}
some
$ y\in \pd D $ \blu{such that} $ K=D\cap H(x_1,y\cdot x_1) $ is a {$(N-1)$-dimensional} face of \reed{$D$.} 
{We} 
\reed{may}
find 
$ u_1,u_2\in K $ and $  {r}\neq 0 $  such that $ u_1-u_2= {r}(x_3-x_2) $. We will use Proposition \ref{pro:fundamental} 
to deduce that  {the problem:} 
\be
\inf\{\rho(x_0-x):x\in D\}\nonumber
\ee
admits more than one  {minimizer.}
\end{proof}

 {
\br
When the set $W$ in \blu{Lemma \ref{lem:nonuniqueness=>}} is  {nonempty},
it means that an edge of the closed unit ball  {associated to} $\rho$ is
orthogonal to one of its vertices.
\Er
}

\bt \label{thm:1}
Let $ \rho $ be a norm on \blu{the 2-dimensional Euclidean space} $  {E} $ \blu{such that its closed  unit ball 
is the convex hull of a finite number of points.}
 {The} {p}roblem:
\be
\inf \{\rho(x_0-x)+\rho^*(x):x\in E\}\nonumber
\ee
admits a unique solution for all $ x_0\in E $  if and only if the following
set:
\begin{gather}
W=\left\{(x_1,x_2,x_3)\in\textnormal{extr}(B _{\rho}(0,1))\times\textnormal{extr}(B _{\rho}(0,1))\times\textnormal{extr}(B _{\rho}(0,1)),\right.\nonumber\\
\left[x_2,x_3]\p  {\pd \bar B_\rho}(0,1),\; x_2\neq x_3,\; x_1\cdot(x_2-x_3)=0\right\}\nonumber
\end{gather}
is empty.
\et

\begin{proof}
If $ W $ is  {nonempty}, by  \blu{{L}emma \ref{lem:nonuniqueness=>}}, we have  {more than one} minimizers.
Suppose  {now that} we have more than one minimizer{s} and  {let} $ D $  {be the} polar of $ B_\rho(0,1) $.
By  {Proposition} \ref{pro:fundamental}
we can find \blu{distinct points} $ u_1, u_2\in \pd \bar B_\rho(0,1) $ and $ w_1,w_2\in \pd D $ 
\blu{such that}
$ [u_1,u_2]\p \pd B_\rho(0,1) $,
$ [w_1,w_2]\p \pd D(0,1) $
and 
$ w_1-w_2 $ is parallel to $ u_1-u_2 $.
Using the characterization of the polar of a set and the fact that the underlying space is of dimension  {two}, we can find an extreme point $ x_1\in \bar B_\rho(0,1) $
that is orthogonal to 
$ w_1-w_2 $.
We will use again the fact that
the underlying space is of dimension  {two,} to  find two distinct extreme points $ x_2,x_3\in \bar B_\rho(0,1) $
\blu{such that} $ x_2-x_3 $ is parallel to $ u_1-u_2 $ \blu{and $[x_1,x_2]\p \pd \bar B_\rho(0,1)$.} 
It holds that $ x_1\cdot(x_2-x_3)=0 $ 
and  {hence} $ W $ is  {nonempty}.
\end{proof}

\br  
We conjecture that Theorem \ref{thm:1}
can be generalized to dimensions  {greater than two}.
\Er

\section{Numerical Experiments}
One shows that the  the problem:
\be
\min\{J {_d}(u_0-u)+J {_d}^*(u): u\in X\}\nonumber
\ee
admits a unique minimizer  {for all $u$ in $X$} if and only if for all $g_0\in Y$ 
the operator  {div} is constant on the set: 
\be
S(g_0)= {\text{argmin}}\{J( {\text{div}}(h-g_0)): h\in Y,\; \|h\|_{\infty}\leq 1\}.\nonumber
\ee
\blu{We start by picking $g_0$ outside the set $\{g\in Y:\|g\|_\infty\leq 1\}$ but close to it and then choosing randomly a number of points to initialize a projected subgradient algorithm. Next, we check if the \pur{iterates generated by} the algorithm \pur{converge} to the same solution or not}.
{Table} \ref{tab:1} show{s} the results  {for} the  \blu{\texttt{C++} implementation of the algorithm} {with} $ N=16 $  {and} the choice of 
$ 500 $ \blu{initialization} points and 
$ 200000$ iterations of the  \blu{algorithm}.

\begin{table}[h]
\centering
\caption{Diameters of solution sets obtained using the distance induced by the norm $\|\cdot\|_X$.}
\label{tab:1}
\begin{tabular}{|c|c|c|c|c|c|}
\hline
Experiment Number&1&2&3&4&5\\
\hline
Diameter&
0.00123793 &
0.00122184 &
0.00150921 &
0.00126712 &
0.00139588 
\\
\hhline{|=|=|=|=|=|=|}
Experiment Number&6&7&8&9&10\\
\hline
Diameter&
0.00155337&
0.00131911&
0.00135451&
0.00119512&
0.00100439 
\\
\hline
\end{tabular}
\end{table}


\section*{Acknowledgements}
The work presented in this paper was initiated while both the authors were postdoctoral fellows at the Institute for Mathematics and its Applications (IMA) during the IMA's annual program on \textit{``Control Theory and its Applications"}.

\bibliographystyle{abbrv}

\end{document}